\documentclass{article} 
\usepackage{amsmath}
\usepackage{amssymb}
\usepackage{amsthm}
\usepackage{graphicx,psfrag}
\usepackage{cite}

\newtheorem{definition}{Definition}
\newtheorem{remark}{Remark}
\newtheorem{example}{Example}

\newtheorem{theorem}{Theorem}

\newtheorem{proposition}{Proposition}

\title{\LARGE \bf
The Statistical Mechanics of Fluctuation-Dissipation \\ and Measurement Back Action\footnote{H. Sandberg is supported by the Hans Werth\'en foundation and a postdoctoral grant from
the Swedish Research Council. J.-C. Delvenne holds an FNRS fellowship (Belgian Fund
for Scientific Research).  This paper presents research results of
the Belgian Programme on Interuniversity Attraction Poles, initiated
by the Belgian Federal Science Policy Office. The scientific
responsibility rests with its authors.}}

\author{Henrik Sandberg$^1$, Jean-Charles Delvenne$^{1,2}$, and John C. Doyle$^1$ \\[0.3cm]
$^1$California Institute of Technology\\ Control and Dynamical Systems\\ M/C 107-81, Pasadena, CA 91125, USA \\[0.3cm]
$^2$Universit\'{e} catholique de Louvain\\ 4, avenue Lema\^{i}tre \\ B-1348, Louvain-la-Neuve, Belgium}

\begin{document}
\maketitle
\thispagestyle{empty}
\pagestyle{empty}

\begin{abstract}
In this paper, we take a control-theoretic approach to answering
some standard questions in statistical mechanics. A central
problem is the relation between systems which appear
macroscopically dissipative but are microscopically lossless. We
show that a linear macroscopic system is dissipative if and only
if it can be approximated by a linear lossless microscopic
system, over arbitrarily long time intervals. As a by-product,
we obtain mechanisms explaining Johnson-Nyquist noise as initial
uncertainty in the lossless state as well as measurement back
action and a trade off between process and measurement noise.
\end{abstract}

\section{Introduction}

The derivation of thermodynamics as a theory of large systems
which are microscopically governed by fundamental laws of
physics (Newton's laws or quantum physics) has a large
literature and tremendous progress for over a century within the
field of statistical physics. See for instance \cite{wannier}
for a physicist's account of statistical mechanics.
Nevertheless, from a control theorist's perspective, there are
inadequacies in the existing treatment both with the level of
mathematical rigor, and the applicability to
far-from-equilibrium systems, particularly when subject to
complex regulatory mechanisms. Substantial work has already been
done in formulating various results of classical thermodynamics
in a more mathematical framework (e.g.
\cite{Brockett+78,Mitter+05,Haddad+05,barahona+02,Bernstein+02}
is a small sample), but statistical mechanics has received much
less comparable attention. This paper focuses on simple problems
in statistical mechanics in which the issue of rigor can be
pursued, but aims also to set the stage for broader
applicability.

In particular, we construct a simple and clear control-theoretic
modeling framework in which the only assumptions on the nature
of the physical systems are conservation of energy and causality
and all systems are of finite dimension and act on finite time
horizons. We construct high-order lossless systems that
approximate low-order dissipative systems in a systematic
manner, and  prove that a linear model is dissipative if and
only if it is arbitrarily well approximated by lossless causal
linear systems over an arbitrary long time horizon. We show how
the error between the systems depend on the number of states in
the approximation and the length of the time horizon. Since
human experience is based on a finite window of space and time,
we argue that no human can directly distinguish between a
low-order macroscopic dissipative system and its high-order
lossless approximation.

The lossless systems studied here are consistent with classical
physics, since they conserve energy, are causal, and are time
reversible. Uncertainty in their initial state gives a simple
explanation of the \emph{Johnson-Nyquist noise} that can be
observed at a macroscopic level. We also derive some well-known
results from statistical mechanics, including the
\emph{fluctuation-dissipation theorem}. As a further
application, we study the implications of these results for an
idealized measurement device, and exhibit a back-action effect,
that there is no precise measurement without perturbation on the
measured system, that arises naturally in a purely classical
setting.

We hope this paper is a step towards building a framework for
understanding fundamental limitations in control and estimation
that arise due to the physical implementation of measurement and
actuation devices. We defer many important and difficult issues
here such as how to actually model measurement devices
realistically. It is also clear that this framework would
benefit from a behavioral setting \cite{polderman+97}. However,
for the points we make with this paper, a conventional
input-output setting with only regular interconnections is
sufficient.  Aficionados will easily see the generalizations,
the details of which might be an obstacle to readability for
others. Perhaps the most glaring unresolved issue is how to best
motivate the introduction of stochastics. In conventional
statistical mechanics, a stochastic framework is taken for
granted, whereas we aim to explain if and when stochastics arise
naturally, and in this we are only partially successful.

The organization of the paper is as follows: In
Section~\ref{sec:micro}, we define the class of linear
lossless/causal systems. In
Section~\ref{sec:physapprox_statical}, we derive lossless/causal
approximations of memoryless dissipative systems and obtain
Johnson-Nyquist noise. In Sections~\ref{sec:interconnect} and
\ref{sec:backaction}, we discuss interconnections of systems and
introduce an idealized measurement device with back action.
Finally, in Section~\ref{sec:memory} we generalize the procedure
from Section~\ref{sec:physapprox_statical} to a class of linear
dissipative systems with memory, and in Section~\ref{sec:fluct}
obtain the fluctuation-dissipation theorem.

\section{Lossless/Causal Linear Systems}
\label{sec:micro}
In this paper, we consider linear systems in the form
\begin{equation}
\begin{aligned}
  \dot x(t) & = J x(t) + Bu(t), \quad x(t) \in \mathbb{R}^n, \\
  y(t) & = B^T x(t),
\end{aligned}
\label{eq:physlinear}
\end{equation}
where $J=-J^T$ and $(J,B)$ is controllable. It is assumed that the input $u(t)$ and the output $y(t)$ are scalars.
We define
the internal energy of (\ref{eq:physlinear}) as
\begin{equation*}
  U(x(t)) \triangleq \frac{1}{2} x(t)^Tx(t).
\end{equation*}
We argue these systems have desirable ``physical'' properties. These properties are losslessness and causality.

\emph{Lossless} \cite{willems72A,willems72B} means that the internal energy satisfies
\begin{equation*}
  \frac{dU(x(t))}{dt} = x(t)^T\dot x(t)= y(t)^Tu(t) \triangleq w(t),
\end{equation*}
where  $w(t)$ is the \emph{work rate} on the system. If there is no work done on the system, $w(t)=0$, then the
internal energy $U(t)$ is constant and conserved. If there is work done on the system, $w(t)>0$, the internal energy increases.
The work, however, can be extracted again, $w(t)<0$, since the energy is conserved and the system is controllable.
Conservation of energy is
a common assumption on microscopical models in statistical mechanics \cite{wannier}.

\emph{Causal} here means that there is no direct term between the input $u$ and the output $y$. This
means that there is no instantaneous reaction of the system. Also this is a reasonable physical assumption.
\begin{definition}
\label{def:physical}
  Systems (\ref{eq:physlinear}) that satisfy the above assumptions are simply called \emph{lossless/causal systems}.
\end{definition}
Later we will seek approximations of dissipative systems in the class of lossless/causal systems.

The lossless/causal systems are rather abstract but have properties that we argue are reasonable from
a physical point of view, as illustrated by the following example.

\vspace{0.1cm}
\begin{example}
Consider the inductor-capacitor circuit in
Fig.~\ref{fig:LCcircuit}. Let the input $u$ be the current through the current
source, and the output $y$ the voltage across the current source. Then a model is given by
\begin{align*}
\dot x &=
\begin{pmatrix}
0 & -1/\sqrt{C_1 L_1} & 0 \\
1/\sqrt{C_1 L_1}  & 0 &  -1/\sqrt{L_1 C_2} \\
0 & 1/\sqrt{ L_1 C_2} & 0
\end{pmatrix}
x +
\begin{pmatrix}
1/\sqrt C_1 \\ 0 \\ 0
\end{pmatrix}
u, \\ \\
y & =
\begin{pmatrix} 1/\sqrt C_1 & 0 & 0 \end{pmatrix}
x, \quad x^T= \begin{pmatrix} \sqrt C_1 v_1 & \sqrt L_1 i_1 & \sqrt C_2 v_2 \end{pmatrix} \\ \\
U &= \frac{1}{2}x^Tx = \frac{1}{2}(C_1v_1^2 + L_1 i_1^2 + C_2 v_2^2), \quad w = yu = v_1i,
\end{align*}
and it satisfies Definition~\ref{def:physical}.
\begin{figure}
  \centering
  \includegraphics[width=0.7\hsize]{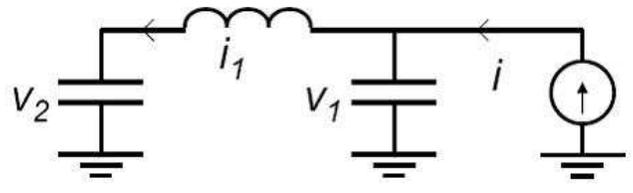}
  \caption{Inductor-capacitor circuit.}
  \label{fig:LCcircuit}
\end{figure}
\end{example}

\section{Lossless/Causal Approximations of Dissipative Memoryless Systems}
\label{sec:physapprox_statical}
In this section, we see how dissipative models, models where energy disappears, can be approximated by
the lossless/causal models. We start with simple memoryless models,
which give rise to heat baths and Johnson-Nyquist noise.
\subsection{Dissipative memoryless systems}
\label{sec:staticio}
Many times macroscopic systems, such as resistors,
can  be modeled approximately by simple input-output relations
\begin{equation}
  y(t) = k u(t),
\label{eq:diss}
\end{equation}
where $k$ is a scalar. If $k>0$, the system is dissipative since we can never extract any work. This is because
the work rate is always
positive
\begin{equation*}
  w(t) = y(t)u(t)= ku(t)^2\geq 0,
\end{equation*}
for all $t$ and $u$. Hence, (\ref{eq:diss}) is neither lossless nor causal.
Next, we show how we can approximate (\ref{eq:diss}) arbitrarily well with a lossless/causal system
over \emph{finite}, but arbitrarily long, time horizons.

First, choose a time interval of interest, $[0,\tau]$, and rewrite (\ref{eq:diss}) using a convolution integral
\begin{equation}
  y(t) = \int_{0}^{\tau} k\delta(t-s) u(s) ds,
\label{eq:resmod}
\end{equation}
when $u$ is at least continuous and has compact support on $[0,\tau]$, and $\delta$ is the Dirac distribution.
Let us call $\tau$ the \emph{recurrence time} of the model.
The recurrence time interval
contains all the time instants where we perform experiments on the model,
and can be very long. Over this time interval,
the system is equally well modeled by the impulse response
\begin{equation*}
  \kappa(t) = \sum_{l=-\infty}^{\infty} k\delta(t-l 2\tau)
\end{equation*}
which is a $2\tau$-periodic distribution. $\kappa(t)$ can be expanded in a Fourier series with
convergence in the sense of distributions:
\begin{equation}
  \kappa(t) \sim \frac{k}{2\tau}+\sum_{l=1}^{\infty} \frac{k}{\tau} \cos l\omega_0 t,
\label{eq:kdistr}
\end{equation}
where $\omega_0=\pi/\tau$. Define the truncated Fourier series by
\begin{equation*}
    \kappa_N(t)\triangleq \frac{k}{2\tau}+\sum_{l=1}^{N} \frac{k}{\tau} \cos l\omega_0 t.
\end{equation*}
We can split $\kappa_N(t)$ into its causal and anti-causal parts:
\begin{align*}
  \kappa_N(t) & \triangleq \kappa_N^c(t) + \kappa_N^{ac}(t) \\
  \kappa_N^c(t) & = 0, \quad t<0 \\
  \kappa_N^{ac}(t) & = 0, \quad t\geq 0.
\end{align*}
We can realize the causal part $\kappa_N^c(t)$ as the impulse response of a lossless/causal system
of order $2N+1$ with the matrices
\begin{equation}
\begin{aligned}
J_N & = \begin{bmatrix} 0 & \Omega_N &  0 \\ -\Omega_N^T & 0 & 0 \\ 0 &0 &0\end{bmatrix},
\,\, \Omega_N = \text{diag}\{\omega_0,2\omega_0,\ldots, N\omega_0\},\\
  C_N & = \sqrt{\frac{k}{\tau}} \begin{pmatrix} 1 & \ldots & 1 & 0 & \ldots & 0 & \dfrac{1}{\sqrt{2}}
\end{pmatrix} \\
  B_N & = C_N^T.
\end{aligned}
\label{eq:harmonics}
\end{equation}
We can realize $\kappa_N^{ac}(t)$ with a similar system by reversing time. That the series (\ref{eq:kdistr})
converges in the sense
of distributions means that for all smooth $u$ of compact support on $[0,\tau]$ we have that
\begin{equation*}
  ku(t) = \lim_{N\rightarrow \infty} \left( \int_0^{\tau} \kappa_N^{ac}(t-s)u(s)ds + \int_0^{\tau}\kappa_N^c(t-s)u(s)ds
\right).
\end{equation*}
A closer study of the two integrals reveals that
\begin{align*}
   \lim_{N\rightarrow \infty} \int_0^{\tau} \kappa_N^{ac}(t-s)u(s)ds & = \frac{1}{2} k u(t+), \\
   \lim_{N\rightarrow \infty} \int_0^{\tau}\kappa_N^c(t-s)u(s)ds & = \frac{1}{2} k u(t-),
\end{align*}
because of the anti-causal/causal decomposition. Hence, since $u$ is continuous, we can model $y(t)=ku(t)$
with only the causal part if we normalize the causal part with a factor two.

We identify the lossless/causal approximation of (\ref{eq:diss}) with a linear operator $K_N:\mathcal{C}^2(0,\tau)\rightarrow
\mathcal{C}^2(0,\tau)$:
\begin{equation*}
  y_N(t)=K_Nu(t): \quad y_N(t)= \int_0^t 2\kappa_N^{c}(t-s)u(s)ds.
\end{equation*}
It is realized by the triple $(J_N,\sqrt{2}B_N,\sqrt{2}C_N)$.
We can bound the approximation error as seen in the following proposition.
\begin{proposition}
\label{prop:errbound}
Assume that $u\in\mathcal C^2(0,\tau)$ and $u(0)=0$. Let
$y(t)=ku(t)$, $k>0$,  and $y_N(t)=K_Nu(t)$. Then
\begin{equation*}
  |y(t)-y_N(t)| \leq \frac{2k\tau}{\pi^2 N}\left( |\dot u(t)|+|\dot u(0)|+\|\ddot u \|_{L_1[0,t]}\right),
\end{equation*}
for $t$ in $[0,\tau]$.
\end{proposition}
\begin{proof}
We have that
\begin{equation*}
  y(t)-y_N(t) = \sum_{l=N+1}^{\infty}\frac{2k}{\tau} \int_0^{t} \cos l\omega_0(t-s)u(s)ds, \,\, t\in[0,\tau].
\end{equation*}
We have changed the order of summation and integration because this is how the value of the series is defined
in distribution sense. We proceed by using repeated integration by parts on each term in the series. We have
\begin{align*}
  & \int_0^t \cos l\omega_0(t-s)u(s)ds = \int_0^t \frac{\sin l\omega_0(t-s)}{l\omega_0}\dot u(s)ds \\
& = \frac{1}{l^2\omega_0^2}\dot u(t)-\frac{\cos l\omega_0 t}{l^2\omega_0^2}\dot u(0) -\int_0^t \frac{\cos l\omega_0(t-s)}{l^2\omega_0^2}\ddot u(s)ds.
\end{align*}
Hence, we have the bound
\begin{equation*}
  |y(t)-y_N(t)|
\leq \frac{2k}{\tau} \sum_{l=N+1}^{\infty} \frac{1}{l^2 \omega_0^2} \left( |\dot u(t)| +|\dot u(0)|+ \int_0^t |\ddot u(s)| ds \right).
\end{equation*}
Since $\sum_{l=N+1}^{\infty}1/l^2 \leq 1/N$, we can establish the bound in the proposition.
\end{proof}

The proposition shows that by choosing $N$ sufficiently large,
we can approximate the memoryless model (\ref{eq:diss}) as well
as we like with a lossless/causal system, if inputs are smooth.
It is a reasonable assumption that inputs, such as voltages, are smooth since
we usually cannot change them arbitrarily fast due to physical limitations. Physically, we can think of $2N+1$ as the
number of degrees of freedom in a resistor. This is usually a number with the size of Avogadro's number, $N\approx 10^{23}$.
Then the recurrence time $\tau$ can be very large without a significant error.
This explains how the dissipative model
(\ref{eq:diss}) is consistent with a physics based on energy conserving systems.

\subsection{Initial conditions in $K_N$}
The general solution to the lossless/causal approximation $K_N$ is
\begin{equation}
  y_N(t) = \sqrt{2}B_N^Te^{J_Nt}x(0) + \int_{0}^t 2\kappa_N^c(t-s)u(s)ds,
\label{eq:gensol}
\end{equation}
where $J_N$ and $B_N$ are defined in (\ref{eq:harmonics}), and $x(0)$ is the initial state.
It is the second part of the solution that approximates $ku(t)$. The first part,
the homogeneous solution, is not desired in the approximation, but is always present for a
linear dynamical system. Next, we study the
influence of this term.

Proposition~\ref{prop:errbound} suggests that we will need a system of incredibly high order to approximate
the dissipative system (\ref{eq:diss}) on a reasonably long time horizon. When dealing with systems of such
extremely high dimensions, it is reasonable to assume that the exact initial state $x(0)$ is not known. Therefore,
we will take a statistical approach to study its influence.

We have that
\begin{equation*}
  \mathbf{E}y_N(t) = \sqrt{2}B_N^Te^{J_Nt}\mathbf{E}x(0) + \int_{0}^t 2\kappa_{N}^c(t-s)u(s)ds,
\end{equation*}
if the input $u$ is deterministic and known.
The covariance function for $y_N(t)$ is then
\begin{equation}
  R_{y_N}(s,t)\triangleq \mathbf{E}[y_N(t)-\mathbf E y_N(t)][y_N(s)-\mathbf E y_N(s)] \\ = 2B_N^Te^{J_Nt}Xe^{-J_Ns}B_N,
  \label{eq:resnoise}
\end{equation}
where $X$ is the covariance of the initial state,
\begin{equation}
X\triangleq \mathbf{E}[x(0)-\mathbf Ex(0)][x(0)-\mathbf Ex(0)]^T.
\label{eq:covariance}
\end{equation}
In Section~\ref{sec:equipartition},
we discuss how it is reasonable to choose $X$. The arguments are information theoretical and physical
in nature. Both arguments result in an equipartition-type statement that result in the concept of temperature.
For now, let us only define the notion of temperature of a lossless/causal system.
\begin{definition}[Temperature]
\label{def:temp}
  A lossless/causal system with deterministic input has temperature $T$ ($T$ is scalar) if
\begin{equation*}
  R_{y}(s,t) =  T\cdot B^Te^{J(t-s)}B.
\end{equation*}
\end{definition}

If $X$ commutes with $J$ and admits $B_N$ as an eigenvector with eigenvalue $T$, (\ref{eq:resnoise}) satisfies Definition~\ref{def:temp}
and we  have (in the sense of distributions)
\begin{equation}
  R_{y_N}(s,t) \rightarrow 2 T k \delta(t-s),\quad t,s\in[0,\tau], \quad N\rightarrow \infty.
\label{eq:limitcovariance}
\end{equation}
A stochastic signal with this property is called \emph{white noise}.

\subsection{Johnson-Nyquist noise}
>From Proposition~\ref{prop:errbound}, (\ref{eq:gensol}), and (\ref{eq:limitcovariance})
we obtain the following proposition.
\begin{proposition}
\label{prop:johnsonnoise}
In the limit when $N\rightarrow \infty$, the lossless/causal system $K_N$, given by (\ref{eq:gensol}),
converges to
\begin{equation}
  y_{\infty}(t) = ku(t) + \sqrt{2Tk} w(t), \quad t\in[0,\tau],
\label{eq:heatbath}
\end{equation}
when it has temperature $T$. The signal $w(t)$ is stochastic white noise of unit intensity. The input
$u(t)$ should satisfy the assumptions of Proposition~\ref{prop:errbound}.
\end{proposition}
\begin{definition}[Heat bath]
\label{def:heatbath}
A system (\ref{eq:heatbath}) is called a \emph{heat bath} of strength $k$, temperature $T$, and
recurrence time $\tau$.
\end{definition}

Hence, in the limit, the uncertainty in the initial state of the microscopic lossless/causal model $K_N$ is transformed into
white noise added to the output of the macroscopic model (\ref{eq:diss}).
This is a generalization of Johnson-Nyquist noise of resistors, see \cite{johnson28,nyquist28}:
It is a fact that careful measurements of the voltage
across a resistor reveal that there is noise that depends on the resistance and temperature. Usually this noise
is modeled by stochastic white noise.
The noise is often explained using methods from statistical
mechanics and circuit theory. See, for example, \cite{wannier}.
Here we obtain exactly the same result using lossless/causal systems and a suitable definition of
temperature.
\begin{remark}
That Proposition~\ref{prop:johnsonnoise} indeed leads to the standard form of the Johnson-Nyquist
noise of a resistor can be seen as follows: We have $v= R i$ from Ohm's law. Assume that $i=0$ and study
the variance of $v(t)$ through a low-pass filter of bandwidth $B$. Then we have, since $|\hat R_w(j\omega)|^2=1$
(white noise),
\begin{equation*}
\mathbf E v(t)^2 = \int_{-B}^B 2TR|\hat R_w(j\omega)|^2 d\omega = 4TRB,
\end{equation*}
which is usually how Johnson-Nyquist noise is presented. Notice that Boltzmann's constant here should be included
in the temperature $T$. It is also interesting to notice that the factor two in the noise intensity $2TR$
in our derivation originates from the causal/anti-causal decomposition in the construction of $K_N$.
A very different argument is used in the derivation in \cite{wannier}.
\end{remark}

\subsection{Equipartition of energy}
\label{sec:equipartition}
In this section, we discuss how the covariance of the initial state $x(0)$ of $K_N$, defined in (\ref{eq:covariance}),
should be chosen. This discussion leads up to the definition of temperature, Definition~\ref{def:temp}.
The first argument
is information theoretical, and the second argument has a more physical flavor. As mentioned in the introduction,
how to properly motivate the introduction of the stochastic element is not easy. Here we just give two arguments whose
consequences are compatible with macroscopic observations, if Johnson-Nyquist noise is modeled by
stochastic white noise. Neither of the arguments is entirely convincing, and we
hope to return to these issues elsewhere.

\subsubsection*{MaxEnt argument}
The first argument is based on the MaxEnt principle, due to Jaynes \cite{Jaynes57A,Jaynes57B}.
This means that we should assign the distribution of $x(0)$ that
maximizes the Shannon entropy of the distribution subject to all known constraints.
The procedure is justified because it leads to the least biased guess.
Assume that the expected internal energy of the initial state is $E$:
\begin{equation*}
  E = \mathbf E \frac{1}{2}x(0)^Tx(0) = \frac{1}{2}\mathbf Ex(0)^T\mathbf Ex(0)+ \frac{1}{2}\text{Tr}X.
\end{equation*}
Maximization of the Shannon entropy subject to this constraint leads to a distribution of $x(0)$
that is Gaussian with mean zero and with covariance matrix
\begin{equation*}
  X = \frac{2E}{2N+1}\cdot I_{2N+1}.
\end{equation*}
If we define the temperature $T$ as $2E/(2N+1)$ and use this $X$ in (\ref{eq:resnoise}),
we see that the covariance function of $y_N$ satisfies
the requested relation in Definition~\ref{def:temp}. This means that the energy is distributed equally
between all degrees of freedom. We have equipartition. The temperature is the expected amount of energy
(up to a factor two) of each degree of freedom. This coincide with the usual notion of temperature in physics.

\subsubsection*{White noise argument}
Assume that the $K_N$ had temperature zero a long time back, i.e., $x(-h)=0$ where $h$ is a large number. We
will be more precise about the size of $h$ later.
We start our experiment at time $t=0$ and wonder what a reasonable
assumption on the initial state
$x(0)$ is. Let us now assume that $K_N$ has been subject to low-intensity white
noise over the time interval $[-h,0]$, say
\begin{equation*}
  \mathbf E u(t)u(s) = \frac{i}{h}\delta(t-s), \quad \mathbf E u(t) =0,
\end{equation*}
where $i$ is an intensity constant. One can say that $K_N$ has been weakly connected to an even larger
heat bath for a long time.

In the end, we want to compute $R_{y_N}$ as defined in (\ref{eq:resnoise}), and it is of interest to compute $X$.
We have
\begin{align*}
  X =\mathbf Ex(0)x(0)^T & = \frac{2i}{h}\int_{-h}^0 e^{-J_Ns}B_NB_N^Te^{J_Ns}ds \\
   & = \frac{2i}{h}\int_{-h}^0 \frac{k}{\tau} \begin{pmatrix} \cos \omega_0 s \\ \cos 2\omega_0 s \\ \vdots \\
\sin \omega_0 s \\ \sin 2\omega_0 s \\ \vdots \\ 1/\sqrt 2 \end{pmatrix}
\begin{pmatrix} \cos \omega_0 s \\ \cos 2\omega_0 s \\ \vdots \\
\sin \omega_0 s \\ \sin 2\omega_0 s \\ \vdots \\ 1/\sqrt 2 \end{pmatrix}^T ds.
\end{align*}
Notice that if $h=2\tau$ we have that
\begin{equation*}
  X= \frac{i k}{\tau}  I_{2N+1},
\end{equation*}
this is the amount of time the white noise needs to excite all the modes equally.
When $h>2\tau$ we can use that
\begin{align*}
  \lim_{h\rightarrow \infty} \frac{1}{h} \int_{-h}^0 \cos k\omega_0s\,\, \cos l\omega_0 s \, ds & =\frac{1}{2}\delta_{k-l} \\
  \lim_{h\rightarrow \infty} \frac{1}{h} \int_{-h}^0 \sin k\omega_0s\,\, \cos l\omega_0 s \, ds & = 0.
\end{align*}
Hence we have that $X \rightarrow (i k/\tau)  I_{2N+1}, \quad h\rightarrow \infty$, and
from (\ref{eq:resnoise}) we have
\begin{equation*}
  R_{y_N}(s,t)= 2B_N^Te^{J_Nt}Xe^{-J_Ns}B_N = \frac{ik}{\tau} 2B_N^Te^{J_N(t-s)}B_N.
\end{equation*}
According to Definition~\ref{def:temp}, the temperature of $K_N$ is $T=ik/\tau$.

\section{Interconnections}
\label{sec:interconnect}
\begin{definition}
The \emph{physical interconnection} of the lossless/causal system $(J_1,B_1,B_1^T)$ to the
lossless/causal system $(J_2,B_2,B_2^T)$ is given by
\begin{equation*}
\begin{aligned}
 \frac{d}{dt} \begin{bmatrix} x_1 \\ x_2 \end{bmatrix} & =
  \begin{bmatrix}
    J_1 & -B_1B_2^T \\
    B_2 B_1^T & J_2
  \end{bmatrix}
  \begin{bmatrix} x_1 \\ x_2 \end{bmatrix} +
\begin{bmatrix}
B_1 \\ 0
\end{bmatrix}u \\
y & = B_1^T x_1.
\end{aligned}
\end{equation*}
\end{definition}

The physical interconnection is still lossless/causal. The interconnection
makes physical sense if one studies interconnections of circuit or mechanical models, for example.
It is also a neutral interconnection, as defined in \cite{willems72A}.
Motivated by this definition, and that we in  Section~\ref{sec:physapprox_statical} showed
that the lossless/causal system $(J_N,\sqrt{2}B_N,\sqrt{2}C_N)$ converges to a heat bath, we make
the following definition.
\begin{definition}
\label{def:physheat}
The \emph{physical interconnection} of the lossless/causal model $(J,B,B^T)$ to a heat bath of strength $k$,
 temperature $T$, and recurrence time $\tau$, is given by
\begin{equation}
\begin{aligned}
  \dot x(t) &= (J-kBB^T)x(t) + Bu(t) - B\sqrt{2k T }w(t)  \\
  y(t) &= B^T x(t),
\end{aligned}
\label{eq:S1+heatbath}
\end{equation}
for $t\in [0,\tau]$, where $w$ is stochastic white noise of unit intensity.
\end{definition}

Notice that even though $(J,B,B^T)$ is lossless,
when connected to the heat bath, (\ref{eq:S1+heatbath}) looks dissipative since the eigenvalues of
$J-kBB^T$ have negative real parts.

\section{Back Action of Linear Measurements}
\label{sec:backaction}
As a simple application of the results in Section~\ref{sec:physapprox_statical}
and the definitions in Section~\ref{sec:interconnect},
consider the problem of measuring the output $y(t)$ of the lossless/causal system $(J,B,B^T)$. For this purpose, we
define an idealized measurement device
\begin{equation}
  y_m(t) = k_m y(t),
\label{eq:ideal}
\end{equation}
where $k_m>0$ is a scalar, and the signal $y_m(t)$ is such that we can read it out perfectly. With such a measurement
device, we can also read out the output $y(t)=y_m(t)/k_m$ perfectly.

Now we construct a slightly less idealized measurement device by replacing (\ref{eq:ideal})
by a lossless/causal approximation of
(\ref{eq:ideal}). This is a more physical device, as argued before. According to Section~\ref{sec:physapprox_statical},
we obtain
\begin{equation}
   y_m(t) = k_m y(t) + \sqrt{2k_mT_m}w(t),
\label{eq:ym}
\end{equation}
in the limit if the initial state of the measurement device is not perfectly known. $T_m$ is the temperature of the
device, and it is essentially a  heat bath. If we make a physical interconnection of $(J,B,B^T)$ to (\ref{eq:ym}),
we obtain
\begin{equation}
\begin{aligned}
  \dot x(t)  & = (J-k_mBB^T)x(t) - B\sqrt{2k_m T_m}w(t),  \\
  \hat y(t) & \triangleq y_m(t)/k_m  = B^Tx(t)+\sqrt{\frac{2T_m}{k_m}}w(t),
\end{aligned}
\label{eq:measuredS}
\end{equation}
using (\ref{eq:ym}) and Definition~\ref{def:physheat}, where $\hat y(t)$ is an estimate of $y(t)$.
Acting on the (\ref{eq:measuredS}) we have
\begin{align*}
  & \text{process noise: } &p(t) &\triangleq \sqrt{2k_m T_m }w(t) \\
  & \text{measurement noise: }& m(t) &\triangleq \sqrt{\frac{2T_m}{k_m}}w(t).
\end{align*}
The measurement device generates process noise and dissipation. This is called \emph{back action} of measurements.
This is a well-known phenomenon in quantum physics. Here we obtain a similar effect based on lossless/causal
approximations and using physical interconnections. Also notice that it holds that
\begin{equation}
  \mathbf E p(t)m(s) = 2T_m \delta(t-s).
\label{eq:tradeoff}
\end{equation}
The cross-covariance between process and measurement noise is
independent of the amplification $k_m$ of the measurement device. For large $k_m$, we get a good estimate
of $y$, but on the other hand, the process noise gets large. Hence, there is a trade off.
It is only the temperature $T_m$ of the measurement device that controls trade off in (\ref{eq:tradeoff}).

\section{Lossless/Causal Approximations of Dissipative Systems with Memory}
\label{sec:memory}
In this section, we generalize the procedure from Section~\ref{sec:physapprox_statical}
to dissipative systems that have memory. We consider strictly stable linear causal systems $G$
with impulse response $g$. Their input-output relation is given by
\begin{equation}
 y(t) = \int_0^t g(t-s)u(s)ds.
\label{eq:Gio}
\end{equation}
The system (\ref{eq:Gio}) is dissipative with respect to the work rate $w(t)=y(t)u(t)$ if
\begin{equation*}
   \int_0^T y(t)u(t)dt \geq 0,
\end{equation*}
for all $T\geq 0$ and admissible $u(t)$. An equivalent condition, see  \cite{slotine+91}, is that the transfer function is positive real
\begin{equation}
  \text{Re}\,\hat g(j\omega) \geq 0 \quad \text{for all} \quad \omega.
\label{eq:PR}
\end{equation}
Here $\hat g(j\omega)$ is the Fourier transform of $g(t)$.

The following theorem shows that the system (\ref{eq:Gio}) is dissipative if and only
if it can be approximated arbitrarily well
by a lossless/causal system over any finite time horizon $[0,\tau]$.
\begin{theorem}
\label{thm:diss2}
Assume that $G$ is a linear (causal) system with impulse response $g$, such that $g\in L_1 \cap L_2(0,\infty)$ and
$g'\in L_1(0,\infty)$.
Then $G$ is dissipative if and only if for all $\epsilon > 0 $ and $\tau > 0$ there is a lossless/causal
linear system $G_{\tau}$ with impulse response $g_{\tau}$ such that
\begin{equation}
  \|g-g_{\tau}\|_{L_2[0,\tau]} \leq \epsilon.
\label{eq:L2approx}
\end{equation}
\vspace{0.1cm}
\end{theorem}
\begin{proof}
See appendix~\ref{proofofdiss2}.
\end{proof}
Notice that Theorem~\ref{thm:diss2} shows that a large class of dissipative systems
(macroscopic systems) can be approximated by the lossless/causal systems we introduced in Section~\ref{sec:micro}.

\section{Fluctuation-dissipation Theorem}
\label{sec:fluct}
If a lossless/causal system satisfies Definition~\ref{def:temp}, then by definition we have
\begin{equation*}
  R_{y}(s,t) = T \cdot B^Te^{J(t-s)}B.
\end{equation*}
This can be said to be the \emph{fluctuation} of the system.
The response of the lossless/causal system to an impulse $u(t)=\delta(t)$ is
\begin{equation*}
 B^Te^{Jt}B.
\end{equation*}
If the lossless/causal system approximates a dissipative system over $[0,\tau]$, see Theorem~\ref{thm:diss2},
then the impulse response decays over this time interval. This represents the \emph{dissipation} of the system.
The expressions of the fluctuation and dissipation are equal up to a constant, the temperature. This
is a  property that can be observed in physical systems close to equilibrium (and hence can be linearized).

\section{Conclusions}
In this paper, we defined the class of lossless/causal systems and used them to approximate dissipative systems.
We obtained an if and only if characterization and gave explicit error bounds that depend on the time horizon and the order of the approximations. When
applied to memoryless models, we saw that Nyquist-Johnson noise (macroscopic measurable noise) can be explained
by uncertainty in the initial state of a lossless/causal approximation of very high order. We also saw that using
these techniques, it was relatively easy to obtain a back-action effect of measurements. This gave rise to a
trade off between process and measurement noise.

\subsection*{Acknowledgment}
The authors would like to thank Ben Recht for many helpful suggestions and comments.

\appendix
\section{Proof of Theorem~\ref{thm:diss2}}
\label{proofofdiss2}
We first show the 'if' direction. Assume first the opposite: That there are lossless approximations that satisfy
(\ref{eq:L2approx}) even though $G$
is not dissipative. If $G$ is not dissipative, we can find an input $u(t)$ over a finite interval $[0,T]$
such that
\begin{equation*}
  \int_0^T y(t)u(t)dt = -K_1 < 0,
\end{equation*}
i.e., we extract energy from $G$ even though its initial state is zero. Call $\|u\|_{L_1[0,T]}=K_2$ and
$\|u\|_{L_2[0,T]}=K_3$. For any $\tau>T$ and $\epsilon>0$ we thus have
\begin{equation*}
  \int_0^T (y_{\tau}(t)-y(t))u(t)dt \leq  \epsilon K_2 K_3,
\end{equation*}
by the assumption that lossless approximations $G_{\tau}$ exist and using the Cauchy-Schwarz inequality.
But the lossless approximation satisfies
\begin{equation*}
  \int_0^T y_{\tau}(t)u(t)dt = \frac{1}{2}x_{\tau}(T)^Tx_{\tau}(T),
\end{equation*}
since $x_{\tau}(0)=0$.
Hence,
\begin{equation*}
  -\int_0^T y(t)u(t)dt=K_1 \leq \epsilon K_2K_3 - \frac{1}{2}x_{\tau}(T)^Tx_{\tau}(T)  \leq \epsilon K_2K_3.
\end{equation*}
But since $\epsilon$ can be made arbitrarily small, this leads to a contradiction.

To prove the 'only if' direction we will explicitly construct a $G_{\tau}$ that satisfies (\ref{eq:L2approx}).
We first need to make some definitions. Let
\begin{equation*}
  C \triangleq \frac{2}{\pi}(\|g\|_{L_{\infty}}+\|g'\|_{L_1}),
\end{equation*}
that  is finite when $g,g'\in L_1$. Also define
\begin{equation*}
  \delta(t) \triangleq \int_t^{\infty} |g(s)|ds
\end{equation*}
that is a continuously decreasing function that satisfies $\lim_{t\rightarrow \infty} \delta(t)=0$.
We will need that the recurrence time $\tau$ is such that
\begin{equation}
\delta(\tau)\leq \epsilon^2/(8C).
\label{eq:taufix}
\end{equation}
If the chosen
$\tau$ does not satisfy this relation, we can without loss of generality increase it to
the smallest $\tau$ that satisfies this bound. That this has been done will be assumed in the following.

The model $G_{\tau}$ we construct will be based on a truncated version of the impulse response $g_{N,\tau}(t)$ where
\begin{align*}
  g_{N,\tau}(t) & = \frac{a_0}{2} + \sum_{k=1}^N a_k \cos \frac{k\pi t}{\tau}, \quad t\in[0,\tau], \\
  a_k & = \frac{2}{\tau} \int_0^{\tau} g(t)\cos  \frac{k\pi t}{\tau} dt \\
  \|g_{N,\tau}\|_{L_2[0,\tau]}^2 & = \frac{\tau}{4} a_0^2 + \frac{\tau}{2} \sum_{k=1}^N a_k^2 \leq \frac{\tau}{2}\sum_{k=0}^N a_k^2.
\end{align*}
Assume that $\tau$ is fixed as above. Next pick the smallest $N$ such that
\begin{equation}
  \|g-g_{N,\tau}\|_{L_2[0,\tau]} \leq \frac{\epsilon}{2}.
\label{eq:Nfix}
\end{equation}
Such an $N$ always exist since $g\in L_2$ and the $\cos$-terms are a basis in $L_2[0,\tau]$.

Define
\begin{equation*}
  \hat g_{N,\tau}(j\omega) \triangleq \int_0^{\tau} g(t)e^{-j\omega t}dt,
\end{equation*}
and notice that
\begin{equation*}
  a_k = \frac{2}{\tau} \text{Re}\,\hat g_{N,\tau}\left( j\frac{k\pi}{\tau} \right).
\end{equation*}
We have that
\begin{align*}
  \left|\text{Re}\,\hat g(j\omega)-\text{Re}\, \hat g_{N,\tau}(j\omega) \right| & = \left|\text{Re}\,\int_{\tau}^{\infty}g(t)e^{-j\omega t}dt
\right| \\& \leq
\|g\|_{L_1[\tau,\infty)} = \delta(\tau)\leq \epsilon^2/(8C).
\end{align*}
Since, $\text{Re}\, \hat g(j\omega)\geq 0$
for all $\omega$ by (\ref{eq:PR}), we have
\begin{equation}
  a_k \geq -\frac{\epsilon^2}{4C\tau}.
\label{eq:akbound1}
\end{equation}
We will need a second bound on $a_k$ that bounds
the rate of decay to zero. We have
\begin{align*}
  a_k & = \frac{2}{\tau} \int_0^{\tau} g(t)\cos  \frac{k\pi t}{\tau} dt \\ &= \frac{2}{\tau}\left(\left[g(t)\frac{\tau}{k\pi}\sin
\frac{k\pi t}{\tau} \right]_0^{\tau} - \int_0^{\tau} g'(t)\frac{\tau}{k\pi}\sin  \frac{k\pi t}{\tau}dt\right),
\end{align*}
and thus
\begin{equation}
  |a_k| \leq \frac{C}{k},
\label{eq:akbound2}
\end{equation}
independent of $\tau$. Together, (\ref{eq:akbound1}) and (\ref{eq:akbound2}) give
\begin{equation}
  a_k \geq \max\left\{-\frac{\epsilon^2}{4C\tau}, -\frac{C}{k} \right\} \quad \text{for all}\quad k.
\label{eq:negakbound}
\end{equation}
Next, define
\begin{equation*}
  g_{N,\tau}(t) = g_{N,\tau}^+(t) + g_{N,\tau}^-(t),
\end{equation*}
where $g_{N,\tau}^-(t)$ contains all the terms in $g_{N,\tau}(t)$ with \emph{strictly negative
Fourier coefficients}. Notice that $g_{N,\tau}^+$ can be realized with a linear lossless/causal system. Compare with
(\ref{eq:harmonics}). We can bound the worst-case $L_2$-norm of $g_{N,\tau}^-$. Using (\ref{eq:negakbound}) we have
\begin{align*}
  \|g_{N,\tau}^-\|^2_{L_2[0,\tau]} & \leq \frac{\tau}{2}\sum_{k=0}^N a_k^2
\\ &\leq \sum_{k=0}^{\lfloor \frac{4C^2\tau}{\epsilon^2}\rfloor} \frac{\tau}{2}\frac{\epsilon^4}{16C^2\tau^2}
+ \sum_{k=\lfloor \frac{4C^2\tau}{\epsilon^2}\rfloor + 1}^{\infty} \frac{\tau}{2} \frac{C^2}{k^2}
\\ & \leq \frac{4C^2\tau}{\epsilon^2} \frac{\tau \epsilon^4}{32 C^2 \tau^2}
+ \frac{\epsilon^2}{4C^2\tau}\frac{\tau C^2}{2} = \frac{\epsilon^2}{4},
\end{align*}
independent of how large $N$ is.

A lossless/causal approximation that satisfies the bound (\ref{eq:L2approx})
is now given by $g_{\tau}(t)=g_{N,\tau}^+(t)$, where $\tau$ and $N$ were fixed in (\ref{eq:taufix}) and
(\ref{eq:Nfix}). This is because
\begin{equation*}
 \|g-g_{N,\tau}^+\|_{L_2[0,\tau]} \leq \|g-g_{N,\tau}\|_{L_2[0,\tau]}  +
\|\underbrace{g_{N,\tau}-g_{N,\tau}^+}_{=g_{N,\tau}^-}\|_{L_2[0,\tau]}
\leq \frac{\epsilon}{2}+ \frac{\epsilon}{2} = \epsilon.
\end{equation*}
This concludes the proof.

\end{document}